\documentclass[11pt]{article}
\usepackage{}
\usepackage[paperwidth=7.7in,top=1.3in, bottom=1.4in]{geometry}
\usepackage{amsfonts}
\usepackage{amsthm}
\usepackage{amssymb}
\usepackage{amsmath,amsfonts}
\usepackage{mathrsfs}
\usepackage{cases}
\usepackage{latexsym,bm}
\usepackage{indentfirst}
\usepackage{color}
\usepackage{ifpdf}
\usepackage{graphicx}
\usepackage{psfrag}
\usepackage{dsfont}
\usepackage[
pdfauthor={CSY},
pdftitle={Strong generalized Halin graphs},
pdfstartview=XYZ,
bookmarks=true,
colorlinks=true,
linkcolor=blue,
urlcolor=blue,
citecolor=blue,
bookmarks=true,
linktocpage=true,
hyperindex=true
]{hyperref}

\newcommand{\xiaowuhao}{\fontsize{9pt}{\baselineskip}\selectfont}

\newtheorem{THM}{\textbf{Theorem}}[section]

\newtheorem{LEM}{\textbf{Lemma}}[section]

\newtheorem{COR}{\textbf{ Corollary}}[section]
\newtheorem{PRO}{\textbf{Problem}}
\newtheorem{PROP}{\textbf{Proposition}}[section]
\newcommand{\pf}{\textbf{Proof}.\quad}
\newcommand{\spf}{ \emph{Proof}.\quad}

\newtheorem{CLA}{\textbf{Claim}}[section]
\newcommand{\qqed}{\hfill $\blacksquare$\vspace{1mm}}

\newcommand{\ve}{\varepsilon }

\begin{document}

\title{Characterizing forbidden pairs for  hamiltonian squares  }
\author{ Guantao Chen and Songling Shan
\\{\xiaowuhao  Georgia State University, Atlanta, GA\,30303, USA}}
\date{}
\maketitle
\emph{\textbf{Abstract}.}
The square of a graph is obtained by adding additional edges joining all
pair of vertices of distance two in the original graph.
Particularly, if $C$ is a hamiltonian cycle of a graph $G$, then the
square of $C$ is called a hamiltonian square of $G$.
In this paper, we characterize all possible forbidden pairs,
which implies the containment of a hamiltonian square,
in a 4-connected graph. The connectivity condition is
necessary as, except $K_3$ and $K_4$,   the square of a cycle
is always 4-connected.

{\textbf{Keywords}.} Hamiltonian square; Forbidden pair

\vspace{2mm}

\section{Introduction}
In this paper, we only consider simple and finite graphs.
Let $G$ and $H$ be two  graphs. We use
$G\sqcup H$ to denote the vertex-disjoint union of $G$ and $H$ if $G$
and $H$ are vertex disjoint, use $G\cup H $ to denote the union of $G$
and $H$,
and use $G+H$ to denote the join of $G$ and $H$,
which is the graph on $V(G)\cup V(H)$ with edges
including all edges of $G$ and $H$, and all
edges between $V(G)$ and $V(H)$.
%If $G$ is a bipartite graph with bipartitions $A$ and $B$,
%we denote $G$ by $G[A,B]$.
The notation $\overline{G}$ denotes the complement of
$G$; that is, the graph with vertex set $V(G)$ and
edges between all non-adjacent pairs of vertices in $G$.
The {\it square} of a graph is obtained by adding additional edges joining all
pair of vertices of distance two in the original graph.
Particularly, if $C$ is a  hamiltonian cycle of a graph $G$, then the
square of $C$ is called a {\it hamiltonian square} of $G$.
If $G$ contains a hamiltonian square,
 we then say $G$ has an $H^2$.
The earliest problem on hamiltonian square can be traced back to
a conjecture proposed by P\'osa~\cite{Posa-H2conjecture}.
The conjecture states that {\it any $n$-vertex graph with minimum degree at least $\frac{2n}{3}$
contains a hamiltonian square}.
The complete tripartite graph $K_{t,t,t-1}$ has
minimum degree $2(3t-1)/3-1/3$, but has no $H^2$.
So, if true, the conjecture is best possible. In 1973,
Seymour~\cite{Seymour-kthpower-conjecture} made a more general conjecture,
which says that
{\it any $n$-vertex graph with minimum degree at least $\frac{kn}{k+1}$
contains a $k$th power of a hamiltonian cycle}.
Here, the $k$th power of a graph is obtained by joining every pair
of vertices of distance at most $k$ in the original graph.
P\'osa's conjecture is almost completely solved.
%Jacobson showed that for $\delta(G)\ge 5n/6$\,(unpublished),
%Faudree, Gould, Jacobson, and Schelp
%in their unpublished work improved the minimum degree condition
%to $3n/4$.
In 1994, Fan and H{\"a}ggkvist~\cite{MR1271992} showed
P\'osa's conjecture for $\delta(G)\ge 5n/7$.
Fan and Kierstead~\cite{MR1309357}, in 1996, proved
that  for {\it any $\ve>0$, there is a number $m$, dependent only on $\ve$,
such that if $\delta(G)\ge (2/3+\ve)n+m$,
then $G$ contains the square of a Hamiltonian path between every pair of edges}.
This implies that $G$ then also contains the square of a hamiltonian cycle.
The same authors in 1996~\cite{MR1399673}, showed that if $\delta(G)\ge (2n-1)/3$,
then $G$ contains the square of a hamiltonian path.
For  graphs with large orders,  P\'osa's conjecture was solved
by Koml{\'o}s, S{\'a}rk{\"o}zy,  and  Szemer{\'e}di~\cite{Komlos}
in 1996 using the Regularity Lemma and the Blow-up Lemma.
Using the absorbing method in avoiding using the Regularity Lemma,
Levitt, S{\'a}rk{\"o}zy,  and  Szemer{\'e}di~\cite{Levitt}
in 2010 improved the bound on the orders.
In 2011, Ch{\^a}u, DeBiasio, and Kierstead~\cite{MR2846301}
verified P\'osa's conjecture for $n\ge 200,000,000$.
The work, in investigating
P\'osa's conjecture,   was
trying to find an $H^2$ in graphs with high minimum degrees.
We may ask, what about finding an $H^2$ in other classes
of graphs? One such possible class is the class
of graphs forbidding  some given small graphs.

Given a family $\mathcal{F}=\{F_1,F_2,\cdots, F_k\}$ of graphs,
we say that a graph $G$ is {\it $\mathcal{F}$-free} if $G$ contains
no induced subgraph isomorphic to any of $F_i, i=1,2,\cdots,k$.
Particularly, when $\mathcal{F}=\{F\}$, we simply say that
$G$ is {\it $F$-free}.
If $G$ is $\mathcal{F}$-free, then the graphs in $\mathcal{F}$  are called {\it forbidden subgraphs}.
The use of forbidden subgraphs to obtain classes of graphs possessing special properties
has long been a common graphical technique.
A pair $\{R,S\}$ of connected graphs is called a {\it hamiltonian forbidden pair}
if every 2-connected $\{R,S\}$-free graph is hamiltonian.
%Given a forbidden pair $(R,S)$,
%then $(R,S)$ is called a {\it forbidden pair} for a
%hamiltonian cycle if for any  2-connected graph $G$,
%$G$ is $(R,S)$-free implies that $G$ is hamiltonian.
The characterizations for hamiltonian forbidden pairs
were completely done
\,(for example, see\,\cite{Bedrossian}, \cite{Duffus}, and \cite{forbidden_pair_hamiltonian}).
Research has also been done on characterizing the forbidden pairs
for stronger hamiltonicity properties~\cite{forbidden_pair_hamiltonian}, such as panconectivity\,(a graph $G$ of
order $n$ is said to be panconnected if any two vertices of $G$, say $x$ and $y$, are joined
by paths of all possible lengths $l$ from $dist(x, y)$ to $n - 1$),
pancyclicity\,(an $n$-vertex graph is pancyclic if it contains cycles
of length $l$, for each $3\le l\le n$).
In this paper, we define forbidden pairs for hamiltonian squares\,($H^2$).
A pair of connected graphs $\{R,S\}$  is called an \emph{$H^2$ forbidden pair} if
every 4-connected $\{R,S\}$-free graph
has an $H^2$. Further more, we give a full characterization for all the possible
$H^2$ forbidden pairs.
%As far as we know, no one else has studied such
%forbidden pairs before.

\begin{figure}[!htb]
\label{forbidden subgraphs}
\psfrag{k13}{$K_{1,3}$} \psfrag{k14}{$K_{1,4}$} \psfrag{z1}{$Z_1$}
\psfrag{H}{$H$} \psfrag{H0}{$H_0=H+\overline{K_3}$} \psfrag{G0}{$G_0=K_4+\overline{K_3}$}
\psfrag{+}{$+$}
\begin{center}
  \includegraphics[scale=0.5]{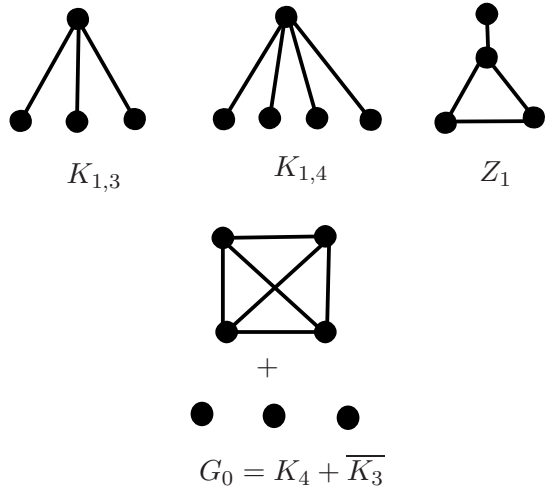}\\
\end{center}
\vspace{-3mm}
  \caption{{\small Small subgraphs}}
  \label{samllgraphs}
\end{figure}

\begin{THM}\label{clawz1}
A pair $\{R,S\}$ of connected graphs with
$R,S\neq P_3$ is an $H^2$ forbidden pair if and only if
$R=K_{1,3}$ and $S=Z_1$, where $Z_1$, as depicted in Figure~\ref{samllgraphs},
 is obtained from
$K_{1,3}$ be adding one edge between two non-adjacent vertices.
%Let $R,S$ be connected graphs\,($R,S\neq P_3$), and let $G$ be a
%4-connected graph. Then $G$ is $(R,S)$-free implies $G$ has
%an $H^2$ if, and only if $R=K_{1,3}$ and $S=Z_1$, where $Z_1$ is
%a triangle attached to one edge as depicted in Figure~\ref{samllgraphs}.
\end{THM}

To force $R=K_{1,3}$ and $S=Z_1$ in Theorem~\ref{clawz1},
a 4-connected 7-vertex graph with no $H^2$  is used in the proof.
Considering graphs with larger order,  we prove
a stronger result.

\begin{THM}\label{k14z1}
A pair  $\{R,S\}$ of connected graphs with
$R,S\neq P_3$  has the property that every 4-connected
$\{R,S\}$-free graph with at least 9 vertices has an $H^2$
if and only if $R\in \{K_{1,3}, K_{1,4}\}$ and $S=Z_1$.
%Let $R,S$ be connected graphs\,($R,S\neq P_3$), and let $G$ be a
%4-connected graph with at least 9 vertices. Then $G$ is $(R,S)$-free implies $G$ has
%an $H^2$ if, and only if $R\in \{K_{1,3}, K_{1,4}\}$ and $S=Z_1$.
\end{THM}

In the study of forbidden pairs for hamiltonian or related properties,
people usually consider pairs $\{K_{1,3}, P_i\}$ for $i\ge 4$.
%In characterizing forbidden pairs for the hamiltonian properties of graphs,
%$\{K_{1,3}, P_i\}$ with some $i\ge 4$ is often one of them.
%So, what about for
%$H^2$?
Except 4 classes of graphs, we show
that all other 4-connected $\{K_{1,3}, P_4\}$-free graphs have an $H^2$,
as given in the theorem below.

\begin{THM}\label{clawp4}
Every
4-connected $\{K_{1,3}, P_4\}$-free graph $G$  has an $H^2$ unless
$G$ is  isomorphic to a graph in one of the following families.
\vspace{-5mm}
\begin{enumerate}
  \item [$($i$)$] $(K_1\sqcup K_3)+(K_m\sqcup K_q)$ with $m+q\ge 4$;
  \item [$($ii$)$] $(K_2\sqcup K_2)+(K_1\sqcup K_m)$ with $m\ge 3$;
  \item [$($iii$)$] $(K_2\sqcup K_3)+(K_1\sqcup K_m)$ with $m\ge 3$;
  \item [$($iv$)$]$(K_3\sqcup K_3)+(K_1\sqcup K_m)$ with $m\ge 3$.
 \end{enumerate}
\end{THM}

It is easy to see that the square of a  cycle  is pancyclic.
%, and contains all
%possible 2-factors.
This is true for any graphs containing an $H^2$. Hence,
partially, we give an answer to a question asked by
Gould  at the 2010 SIAM Discrete Math meeting in Austin, TX.

\begin{PRO}\label{pancyclicity}
Characterize the pairs of forbidden subgraphs that imply a 4-connected graph is pancyclic.
\end{PRO}

It is worth mentioning that
all the known forbidden pairs on Problem~\ref{pancyclicity}
include the claw: $K_{1,3}$\,(see \cite{forbidden_4connected_pancyclicity},
\cite{MR3004480} and \cite{MR2713982}). Hence Theorem~\ref{k14z1}
gives a new forbidden pair for pancyclicity.

\section{Properties of Some Non-hamiltonian Square Graphs}
\begin{figure}[!htb]
\psfrag{km}{$K_{m}$} \psfrag{G2}{$G_{2}\,(m\ge 4)$} \psfrag{G3m=Km+km-1barmle3}{$G_3=K_m+\overline{K_{m-1}}\,(m\ge 4)$}
\psfrag{G4}{$G_4$} \psfrag{G5 in G5(2)}{$G_5\in G_5(2)$} \psfrag{x1}{$x_1$} \psfrag{x2}{$x_2$}\psfrag{x3}{$x_3$}\psfrag{x4}{$x_4$} \psfrag{v4}{$v_4$}
\psfrag{y1}{$y_1$} \psfrag{y2}{$y_2$}\psfrag{y3}{$y_3$}\psfrag{y4}{$y_4$}\psfrag{v6}{$v_6$}
\psfrag{z1}{$z_1$} \psfrag{z2}{$z_2$}\psfrag{z3}{$z_3$}\psfrag{z4}{$z_4$}
\psfrag{G6}{$G_6\,(m\ge 4)$}
\psfrag{+}{$+$}
\psfrag{m}{$m$}
\psfrag{1}{$1$} \psfrag{2}{$2$}
\begin{center}
  \includegraphics[scale=0.4]{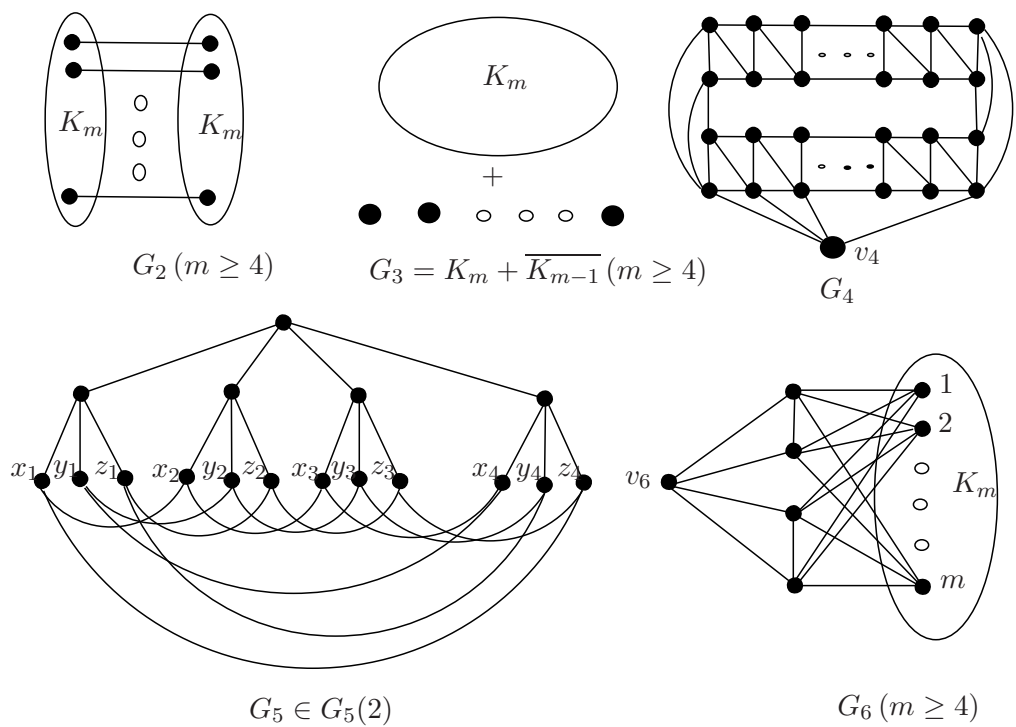}\\
\end{center}
\vspace{-1mm}
  \caption{{\small 4-connected no $H^2$ graphs}}
  \label{forbidden}
\end{figure}

In this section, we examine some  properties of
the graphs depicted  in Figure~\ref{forbidden}.
These graphs will be used in the following
section to characterize  the $H^2$
forbidden pairs. The formal definitions of  these graphs are given below.

\begin{itemize}
\item [$G_1$:] $K_{m,m}$, a complete bipartite  graph with $m$ vertices in each bipartite sets, where $m\ge4$.
  \item [$G_2$:] $K_m\sqcup K_m \cup M$,
  a graph obtained from two vertex-disjoint copies of $K_m$ by adding a perfect matching $M$ between them, where $m\ge4$.
    \item [$G_3$:] $K_m+\overline{K_{m-1}}$, the join of $K_m$ and $\overline{K_{m-1}}$,  where  $m\ge 4$.
  \item [$G_4$:] The graph obtained from the square of a cycle, denoted as $C^2$,
  by  joining a new vertex $v_4$ to four vertices on $C^2$
  such that the four vertices induces $P_3\sqcup K_1$ in the $C^2$.
  \item [$G_5$:]
   Let $T_t$ be a rooted tree of depth $t$\,(the length of a longest path from the root to a leaf is $t$)
  such that all the leaves are at the same depth and all non-leaves have degree 4
  (known as a prefect 4-ary tree). Then  $G_5(t)\,(t\ge 2)$ is  the graph obtained from $T_t$
by connecting the leaves into a cycle in a way such that the girth of the
finally resulted graph is greater than 4. The graph $G_5$ from the family
$G_5(2)$ is depicted in Figure~\ref{forbidden}. $G_5$ is obtained as follows: embed a copy of $T_2$
on the plane, and name the leaves from the left to right, consecutively,
as $x_1,y_1,z_1,x_2,y_2,z_2,\cdots, x_4,y_4,z_4$;
then a cycle $C=x_1x_2x_3x_4y_1\cdots y_4z_1\cdots z_4x_1$
is obtained by joining the corresponding edges.
The construction can be easily generalized to $G_5(t)$
for $t\ge3$. (In $G_5(2)$,  a cycle using  the root vertex
contains three non-leaves and at least two leaves;
and a cycle not using the root vertex
uses at least two non-leaves and 4 leaves.
In any case, it indicates that $G_5(2)$
has girth at least 5. Similarly, $G_5(t)$ has girth at least 5.)
\item [$G_6$: ]$(K_2\sqcup K_2)+(K_m\sqcup K_1)$, where $m\ge 4$.
Denote the isolated vertex in  $K_m\sqcup K_1$ by $v_6$.
\end{itemize}

It is not hard to check that all those graphs are 4-connected.
Furthermore, we have the following fact.

\begin{LEM}\label{no H2}
None of the graphs in Figire~\ref{forbidden} has an $H^2$.
\end{LEM}

\pf Notice that in an $H^2$, the neighborhood  of any vertex induces a $P_4$.
If $G_2$ has an $H^2$, then it must contain one of the edges
connecting the two copies of $K_m$. Let $xy$ be a such edge.
Then the neighbors of $x$ on the $H^2$ consists of $y$
and another three vertices from the copy of $K_m$ containing $x$.
However, those four vertices do not induce a copy of $P_4$,
showing a contradiction. Similarly, neither of the set of neighborhoods of $v_4$
in $G_4$ or of $v_6$ in $G_6$ induces $P_4$. Thus, neither $G_4$
nor  $G_6$ has an $H^2$.
As $G_3=K_m+\overline{K_{m-1}}$,
% has order $2m-1$ and has an independent set of size $m-1$,
any hamiltonian cycle of $G_3$ contains a pair of
vertices from $V(\overline{K_{m-1}})$ such that they have distance 2 on the hamiltonian cycle.
%but are non-adjacent in the original graph $G_3$.
This in turn implies that $G_3$ has no $H^2$.
As an $H^2$ contains triangles, the triangle-free graph $G_5(t)$
has no $H^2$.
\qqed

As the graph $G_2$ will be used more frequently later on,
we discuss its properties in more detail here.

\begin{LEM}\label{G2}
Let $S\notin \{K_3,P_3\}$ be a connected $\{P_4,C_4,K_4\}$-free graph.
If $G_2$ contains $S$ as an induced subgraph,
then $S$ is $Z_1$.
\end{LEM}

\pf Since $V(G)\ne \emptyset$ and $E(G)\ne \emptyset$, $S\notin \{K_1,K_2\}$. Thus $|V(S)|\ge 3$.
Since $S\notin\{K_3,P_3\}$ and any connected
3-vertex subgraph of $G_2$ is  either $K_3$ or $P_3$, we conclude that $|V(S)|\ge 4$.
Furthermore,
as $S$ is $K_4$-free, it contains at most 3 vertices
from one of the copies of $K_m$.
Since $S$ is connected and $\{P_4, C_4\}$-free, if it contains
at least two vertices from one copy of $K_m$, then it
contains at most one vertex from the other copy of
$K_m$. Hence
$S$ contains exactly three vertices from
one copy of $K_m$, and exactly one vertex from
the other.  The connected graph induced on such four
vertices can only be isomorphic to $Z_1$.
\qqed

\section{Proofs of the Main Results}

In this section, we prove Theorem~\ref{clawz1}, Theorem~\ref{k14z1}, and Theorem~\ref{clawp4}.
We first characterize the single forbidden subgraph for 4-connected graphs containing an $H^2$.
As any $P_3$-free graph is complete, we observe that
any 4-connected $P_3$-free graph has an $H^2$. Conversely, we have the following result.

\begin{PROP}\label{p3}
A connected graph $F$ has the property that every 4-connected
$F$-free graph has an $H^2$ if and only if $F=P_3$.
%Suppose that $A$ is a connected graph and $G$ is a 4-connected graph. Then $G$ is $A$-free implies
%$G$ has an $H^2$ only if $A=P_3$.
\end{PROP}

\pf Since $G_1=K_{m,m}$ has no $H^2$, $G_1$ contains $F$ as an induced subgraph.
Hence $F=K_{1,r}$, where $r\ge 2$ or $F$ contains an induced $C_4$. As the graph $G_4$ in
Figure~\ref{forbidden} has no
$H^2$ and is $C_4$-free, we see that $F=K_{1,r}$. The only induced star contained in all the graphs of family $G_2$
is $K_{1,2}$; that is, an induced copy of $P_3$. Hence $F=P_3$.
\qqed

%Let $G$ be a graph.  Recall that a matching  of $G$ is a set of independent edges of $G$.
%We define a mixed  $(1,r)$-matching of $G$ as a set of vertex
%disjoint cliques size between $2$ and $r+1$. Particularly, a mixed $(1,1)$-matching is
%just a matching.

We study the structure of a connected $Z_1$-free graph in the
following theorem, which will help us in knowing
the structure of  a $\{K_{1,r},Z_1\}$-free graph\,($r\ge 3$).

\begin{LEM}\label{z1}
Let $G$ be a connected $Z_1$-free graph.
If there exists a vertex $v\in V(G)$
such that $d(v)\ge 3$ and $v$ is contained in a
triangle, then
$G$ is isomorphic to  a complete multipartite graph $K_{t_1,t_2,\cdots,t_k}$.
\end{LEM}

\pf We use induction on $n=|V(G)|$.
When $n=4$, $G$ is either $K_4$
or the graph obtained from $K_4$ by removing one edge,
so the result holds. Suppose that $n\ge 5$ and
that  Lemma~\ref{z1} holds for graphs with less than
$n$ vertices.
Let $v\in V(G)$ be a vertex such that
$d(v)\ge 3$ and  $v$ is contained in a
triangle.
%We partition $V(G)$ into
%two sets: $N[v]\sqcup\overline{N}[v]$,
Let $N[v]:=N(v)\cup\{v\}$ and  $\overline{N}[v]=V(G)-N[v]$.
Notice that $\overline{N}[v]$ may be
empty. As $G$ is $Z_1$-free,
we know $G[N(v)]$ is $(K_2\sqcup K_1)$-free.
Together with the fact that
$G[N(v)]$ contains an edge,
we then know $G[N(v)]$ is connected.
Before examining the structure of
$G[N(v)]$ further, we
claim the following.

\begin{CLA}\label{remaining}
If $\overline{N}[v]\ne \emptyset$,
then for every $w\in\overline{N}[v]$,
$N(w)=N(v)$ holds.
%$N(v)\subseteq N(w)$ and $\overline{N}[v]$ is
%an independent set in $G$.  This implies
%that $N(v)=N(w)$.
\end{CLA}

\spf  Let $w\in\overline{N}[v]$. We first claim that
if $N(w)\cap N(v)\neq \emptyset$, then $N(v)\subseteq N(w)$.
%We first claim  that
%for any vertex in $\overline{N}[v]$,
%if it is adjacent to one vertex in $N(v)$,
%then it is adjacent to all vertices in $N(v)$.
Suppose not, then
%since $G$ is connected,
there exists
%$%w\in\overline{N}[v]$ such that
%$w$ is adjacent to some vertex in $N(v)$.
%Suppose that $N(v)-N(w)\ne \emptyset$.
 $v'\in N(v)$ such that
$wv'\notin E(G)$. We choose a such $v'$ such that
$w$ is adjacent to a neighbor of $v'$,  say $u'$,
 in $N(v)$.
However, the graph induced on $\{v,v',u',w\}$
is isomorphic to $Z_1$, showing a contradiction.
Hence $N(v)\subseteq N(w)$. The claim is proved.

We then claim that if $N(w)\cap N(v)\neq \emptyset$, then $N(w)\subseteq N(v)$.
%for every $w\in\overline{N}[v]$,
%if $w$ is adjacent to a vertex in $N(v)$, then
%$w$ is not adjacent to any vertex in
%$\overline{N}[v]$.
Otherwise, assume that $w$ is adjacent to
a vertex $w'\in\overline{N}[v]$.
If $w'$ is adjacent to a vertex in
$N(v)$, then we have $N(v)\subseteq N(w)\cap N(w')$
by the earlier assertion.
Let $v'\in N(v)\subseteq N(w)\cap N(w')$.
Then   $\{v,v',w,w'\}$ induces a $Z_1$.
Hence we assume $w'$ is not adjacent to
any vertex in $N(v)$.
Let $v',u'\in N(v)\subseteq N(w)$.
Then $\{u',v',w,w'\}$ induces a $Z_1$.
Thus $w$ is not adjacent to any vertex in
$\overline{N}[v]$.

As $G$ is connected, Claim~\ref{remaining} is then implied by the
above two assertions.
\qed

We now proceed with the proof according to
several cases depending on the structure
of $G[N(v)]$.
Let $|V(G)-N(v)|=t'$ and $G'=G[N(v)]$.
Recall that $G'$ is connected and is $(K_2\sqcup K_1)$-free.

{\indent \it Case 1.}  $G'$ has a vertex with
degree at least 3 in $G'$ and the vertex is contained
in a triangle in $G'$.

By the induction hypothesis, $G'\cong K_{t_1,t_2,\cdots,t_{k-1}}$.
Then we have
$G\cong K_{t_1,t_2,\cdots,t_{k-1},t'}$.

So we suppose that the condition in Case 1 is
not satisfied by $G'$.
Let $u\in V(G')$ be a vertex of maximum degree in $G'$.

{\indent \it Case 2.}  $d_{G'}(u)\le 2$.

Then $G'$ is the union of vertex disjoint paths and cycles.
As $G'$ is connected and is $(K_2\sqcup K_1)$-free,
we know $G'$ is isomorphic to one of the graphs $K_3$, $P_3$, or $C_4$.
In any case,  $G$ is isomorphic to a complete
multipartite graph.

{\indent \it Case 3.}  $d_{G'}(u)\ge  3$.

As $u$ is not on a triangle in $G'$, $N_{G'}(u)$
is an independent set in $G'$. If $N_{G'}[u]=V(G')=N(v)$,
then it is already seen that $G$ is isomorphic to a complete
multiple graph with the size of each parts as $t', 1$,
and $d_{G'}(u)$, respectively. Hence, we assume $N(v)-N_{G'}[u]\ne \emptyset$.
As $G'$ is connected and is $(K_2\sqcup K_1)$-free,
 every vertex in $N(v)-N_{G'}[u]$ is adjacent
to every vertex in $N_{G'}(u)$. Again, by the fact that
$G'$ is $(K_2\sqcup K_1)$-free, we know there is no edge
with the two ends in  $N(v)-N_{G'}[u]$.
Hence, $N(v)-N_{G'}[u]$ is an independent set. Let
$t_1=d_{G'}(u)$ and $t_2=|N(v)-N_{G'}(u)|$.  We see
that $G\cong K_{t_1,t_2,t'}$.

The proof is complete.
\qqed

Additionally, if $G$ is
a $\{Z_1,K_{1,r}\}$-free graph with a vertex of degree at least $r$\,($r\ge 3$),
then $G$ contains a vertex which is contained in a triangle and
is of degree at least 3. Thus by applying Lemma~\ref{remaining}
and by the fact that $G$ is $K_{1,r}$-free, we have the following result.

\begin{COR}\label{k1rz1}
Let $G$ be a connected $\{Z_1,K_{1,r}\}$-free graph
with a vertex of degree at least $r$.
Then
$G$ is isomorphic to  a complete multipartite graph $K_{t_1,t_2,\cdots,t_k}$
such that each $1\le t_i\le r-1$.
\end{COR}

The case of $r=3$ in the above Corollary has been mentioned
in other research papers, for example, in~\cite{forbidden_pair_hamiltonian}.
By Corollary~\ref{k1rz1}, we have the following result.
\begin{COR}\label{connect}
A connected $\{K_{1,r}, Z_1\}$-free graph with a vertex of degree at least $r$
is $(n-r+1)$-connected.
\end{COR}

By Corollary~\ref{k1rz1}, a 4-connected $\{Z_1,K_{1,3}\}$-free graph $G$
is a complete graph missing  at most a matching. By finding
a hamiltonian cycle of $G$ such that non-adjacent pairs of vertices are of
distance at least 3 on the cycle, we can construct  an $H^2$ in $G$.  Hence,
we obtain the result below.

\begin{THM}\label{clawk13}
Every  4-connected $\{Z_1,K_{1,3}\}$-free graph contains an $H^2$.
\end{THM}

For 4-connected $\{Z_1,K_{1,4}\}$-free graphs, we
have a similar result.

\begin{THM}\label{clawk14}
Every  4-connected $\{Z_1,K_{1,4}\}$-free graph
 contains an $H^2$ provided $|V(G)|\ge 9$.
\end{THM}
\pf
Let $n=|V(G)|$. We use induction on $n$ to show the theorem. By Corollary~\ref{k1rz1}, any 4-connected
9-vertex $\{Z_1,K_{1,4}\}$-free graph contains $K_{3,3,3}$
as a spanning subgraph. It is not difficult to verify that
$K_{3,3,3}$ contains an $H^2$. For example,
let $\{x_i,y_i,z_i\}\,(i=1,2,3)$ be the
three vertices in the $i$-th tripartition. Then $x_1x_2x_3y_1y_2y_3z_1z_2z_3x_1$
with the additional edges
gives an $H^2$. So we assume $n\ge 10$.
Let $v\in V(G)$ be a vertex. We consider the graph $G'=G-v$.
Then $G'$ is 6-connected by
Corollary~\ref{connect}.
Additionally,  $G'$ has at least 9 vertices and is $\{Z_1,K_{1,4}\}$-free.  Hence it contains an $H^2$,
say $C_1^2$
by the induction hypothesis. Since $G$
is a multipartite graph with each partition of
size at most 3, there are at most two vertices
on $C_1^2$ which are not adjacent to $v$. Thus,
there are at least 4 consecutive vertices on $C_1^2$
such that each of them is adjacent to $v$.
Let $v_1,v_2,v_3,v_4$ be 4 such consecutive vertices
on $C_1^2$. Then $C_1^2-\{v_2v_3,v_2v_4,v_1v_3\}\cup \{vv_i\,|\,i=1,2,3,4\}$
gives an $H^2$ of $G$.
\qqed

Notice that the order 9 condition in the above theorem is sharp.
The complete tripartite 8-vertex graph $K_{2,3,3}$  is 4-connected and
$\{K_{1,4},Z_1\}$-free, but contains no $H^2$.

Before proving Theorem~\ref{clawz1} and Theorem~\ref{k14z1}, we notice that
if $\{R,S\}$ is a forbidden pair implying  the containment of an $H^2$
in a 4-connected graph, then neither of $R$ or $S$ is a triangle
since an $H^2$ always contains  triangles.

\subsection{Proof of Theorem~\ref{clawz1}}

The sufficiency follows from Theorem~\ref{clawk13}.

Conversely, we will first show that one of $R$ and $S$ must be a claw.
Thus, suppose that $R,S\ne K_{1,3}$. Assume, without loss of generality, that $R$ is an induced subgraph
of $G_1=K_{m,m}$. Then $R=K_{1,r}$, where $r\ge 4$ or $R$ contains an induced $C_4$. We now consider two cases.

\emph{\textbf{Case 1: }} $R=K_{1,r}\,(r\ge 4)$.

The graph $G_4$ has no induced copy of $R$, so it contains an induced copy of $S$.
As $G_4$ is $\{K_4, K_{1,3}\}$-free, we see that
$S$ contains no $K_4$ and no induced $K_{1,3}$.
Also, $R$ is not an induced subgraph of $G_0=K_4+\overline{K_3}$.
So $G_0$ contains $S$ as an induced subgraph.
Since $S\notin \{P_3, K_3\}$ and any  connected 3-vertex
subgraph of  $G_0$ is contained in
$\{P_3, K_3\}$, we conclude that
$S$ has at least 4 vertices.
In $G_0$, any 4
vertices of $G_0$ with at most one vertex in
$\overline{K_3}$ induces a $K_4$;
and any
4
vertices of $G_0$ with three vertices in
$\overline{K_3}$ induces a $K_{1,3}$.
Hence, $S$ contains exactly two vertices
from  the subgraph $K_4$ of $G_0$ and
exactly two vertices from the subgraph
$\overline{K_3}$ of $G_0$, as
$S$ contains no $K_4$, and no induced $K_{1,3}$.
%we conclude that $S$ contains at least two
%vertices from $\overline{K_3}$ by noticing that
%$S$ does not contain   $K_4$.
%As $S$ does not contain   induced $K_{1,3}$,
%it contains at most two vertices from
%the subgraph graph $\overline{K_3}$ of $G_0$.
%Hence, $S$ contains exactly two vertices
%from  $\overline{K_3}$. Again, using that
%fact that  $S$ does not contain   $K_4$ and
%$S$ has order at least 4,
%it contains exactly two vertices from the
%clique $K_4$.
%Putting those argument together, we see that
So $S$ is an
induced $K_4^-$\,($K_4$ with exactly one edge removed).
However, $G_2$ has no  induced $R=K_{1,r}\,(r\ge 4)$,  and no induced $K_4^-$.
We obtain a contradiction.

\emph{\textbf{Case 2: }} $R$ contains an induced $C_4$.

Since $G_4$ has no induced copy of $R$,  it contains an induced copy of $S$.
As $G_4$ is $\{K_4, K_{1,3}\}$-free, we see that
$S$ contains no $K_4$ and no induced $K_{1,3}$.
Also, $R$ is not an induced subgraph of $G_3$.
So $G_3$ contains $S$ as an induced subgraph.
 Since $S$ is connected and $S\notin \{K_1,K_2, K_3,P_3\}$,
 and  any connected  2-vertex, 3-vertex subgraphs of $G_3$
are contained in $ \{K_2, K_3,P_3\}$, we
conclude that $|V(S)|\ge 4$.
In $G_3$, any 4 vertices from $K_m$ or any 3 vertices from
$K_m$ and one vertex from $\overline{K_{m-1}}$
induce a $K_4$; and any
4 vertices in which three from $\overline{K_{m-1}}$
induce a $K_{1,3}$. We conclude
 that $S$ contains exactly two vertices from
$K_m$ and exactly two vertices from $\overline{K_{m-1}}$,
as $S$ contains no $K_4$ and no induced $K_{1,3}$.
%Then $S$ contains at most three vertices from
%the copy $K_m$ of $G_3$  as  $S$ contains
%no $K_4$. Hence $S$ contains
%vertices from both the subgraphs $K_m$ and $\overline{K_{m-1}}$ of $G_3$.
%Notice that any connected subgraph of $G_3$ with one vertex from $\overline{K_{m-1}}$,
%and at least three vertices from $K_m$ contains a $K_4$.
%As $S$ contains no $K_4$, it then implies that $S$ contains
%at least two vertices from $\overline{K_{m-1}}$ and at most two vertices from $K_m$.
%Notice also that any connected  subgraph of $G_3$ with at least three vertex from $\overline{K_{m-1}}$,
%must contain  at least one vertex from $K_m$; and thus contains an induced copy of $K_{1,3}$.
% As $S$ has no induced $K_{1,3}$,  it contains exactly two vertices from $\overline{K_{m-1}}$.
%Recall the argument earlier that $|V(S)|\ge 4$  and $S$ contains
%at most two vertices from $K_m$, we see that $S$ contains exactly two vertices from
%$K_m$ and exactly two vertices from $\overline{K_{m-1}}$.
So  $S$ is an induced $K_4^-$.
However, each graph in   $G_5(t)$ has no $H^2$,  no induced $C_4$,
and  no triangle\,(so no $K_4^-$). This gives a contradiction.

Thus, one of $R$ and $S$ must be a claw. We assume, without loss of generality, that
$R=K_{1,3}$. As $R$ is $K_{1,3}$, $S$ in an induced subgraph of  $G_2$,
$G_4$,  and $G_6$, as none of them  contains induced claws. Note that $G_4$ is $\{C_4,K_4\}$-free,
and $G_6$ is $P_4$-free, so $S$ is $\{P_4,C_4,K_4\}$-free.
Applying  Lemma~\ref{G2}, we see $S$ is $Z_1$.
\qqed

\subsection{Proof of Theorem~\ref{k14z1}}

The sufficiency follows from Theorem~\ref{clawk14}.

Conversely, we will first show that one of $R$ and $S$ must be $Z_1$.
Thus, suppose that $R,S \ne Z_1$.  Assume, without loss of generality, that $R$ is an induced subgraph
of $G_1=K_{m,m}$. Then $R=K_{1,r}$, where $r\ge 3$ or $R$ contains an induced $C_4$. We now consider two cases.

\emph{\textbf{Case 1: }} $R=K_{1,r}\,(r\ge 3)$.

Then $R$ is not an induced subgraph of $G_2$. So $G_2$ contains
$S$ as an induced subgraph.
Both $G_4$ and  $G_6$ contains an induced copy of $S$
since neither of them contains  an induced copy of $R$.
Since $G_4$ is $\{C_4, K_4\}$-free and $G_6$ is $P_4$-free, we see that
$S$ is $\{P_4,C_4,K_4\}$-free.
Applying Lemma~\ref{G2}, we have $S=Z_1$.

\emph{\textbf{Case 2: }} $R$ contains an induced $C_4$.

The graph $G_4$ has no induced copy of $R$, so it contains an induced copy of $S$.
As $G_4$ is $\{K_4, K_{1,3}\}$-free, we see that
$S$ contains no $K_4$ and no induced $K_{1,3}$.
Also, $R$ is not an induced subgraph of $G_3$.
So $G_3$ contains $S$ as an induced subgraph.
Since $S$ is connected and
 $S\notin \{K_1,K_2, K_3,P_3\}$,
and any connected  2-vertex, 3-vertex subgraphs of $G_3$
are contained in $ \{K_2, K_3,P_3\}$, we
conclude that $|V(S)|\ge 4$.
In $G_3$, any 4 vertices from $K_m$ or any 3 vertices from
$K_m$ and one vertex from $\overline{K_{m-1}}$
induce a $K_4$; and any
4 vertices in which three from $\overline{K_{m-1}}$
induce a $K_{1,3}$. We conclude
 that $S$ contains exactly two vertices from
$K_m$ and exactly two vertices from $\overline{K_{m-1}}$,
as $S$ contains no $K_4$ and no induced $K_{1,3}$.
So  $S$ is an induced $K_4^-$.
However, $G_2$ has no  induced $R=K_{1,r}\,(r\ge 3)$  and no induced $K_4^-$.
We obtain a contradiction.

Thus one of $R$ and $S$ must be $Z_1$. Assume, without loss of generality, that
$S=Z_1$. As $G_1=K_{m,m}$ contains no $Z_1$, $G_1$ contains
an induced copy of $R$. Hence $R=K_{1,r}$, where $r\ge 3$ or $R$ contains an induced $C_4$.
Since each graph in $G_5(t)\,(t\ge 2)$ is  $C_4$-free, and the only possible stars
in it are $K_{1,r}$ for $r\le 4$, we see that $R=K_{1,r}$ for $r=3,4$.
\qqed

\subsection{Proof of Theorem~\ref{clawp4}}

We now prove Theorem~\ref{clawp4}.
Let $P$ be a path. We use $P^2$ to denote the square of $P$. In
omitting the edges joining distance
2 vertices on the path, we will use
the same notation to denote the  square of the path.  Similar notation for
the square of a cycle.
Let $P_1^2=v_1v_2\cdots v_{s-1}v_s$ and $P_2^2=u_1u_2\cdots u_{t-1}u_t$
be two path squares.   We denote by $P_1^2P_2^2$  as
the concatenation of $P_1^2$ and $P_2^2$
by adding edges $u_1v_s, u_1v_{s-1}$ and $u_2v_s$,
where $u_1v_{s-1}$ exists  only if $s\ge 2$ and $u_2v_s$
exists only if $t\ge 2$.  Also, the notations $v_1P_1^2, P_1^2v_s$, or
$v_1P_1^2v_s$ may be used for specifying the end vertices
of $P_1^2$.

We may assume that $G$ is not complete. Let $S$ be a minimum vertex-cut of $G$.
Let $G_i=(V_i,E_i)\,(i=1,2,\cdots,k)$ be all the components of $G-S$. Since $G$
is 4-connected, $|S|\ge4$.
As $S$ is a minimum vertex-cut, we have the following claim.

\emph{\textbf{Claim 1:}} For every vertex $v\in S$, $N(v)\cap V_i\ne \emptyset$, for all
$i=1,2,\cdots,k$.

Since $G$ is claw-free, from Claim 1 we get Claim 2 below.

\emph{\textbf{Claim 2:}} $k=2$; that is,  $G-S$ has exactly two components.

Also, by the fact that $G$ is $P_4$-free, we conclude the following claim.

\emph{\textbf{Claim 3:}} For each $v\in S$, $N_{G_i}(v)=V_i$ for $i=1,2$.

As  $E(V_1,V_2)=\emptyset$, $G$ is claw-free, and by  Claim 3, we obtain Claim 4
as follows.

\emph{\textbf{Claim 4:}} $G_i$ is a complete subgraph of $G$ for $i=1,2$.

We will use induction on $n=|V(G)|$ in some cases of the proof.
The smallest 4-connected $\{K_{1,3}, P_4\}$-free graph is $K_5$, it contains an $H^2$.
So we suppose $n\ge 6$ and suppose that the theorem holds for the described graphs
of smaller orders. Let $P_i^2$ be a hamiltonian path square of $G_i\,(i=1,2)$.

If $G[S]$ is 4-connected and is not isomorphic to any graphs in
the exception families, then by the induction hypothesis,
$G[S]$ contains an $H^2$, say $C_s^2$,
which contains at least 4 vertices by the assumption that $G[S]$
is 4-connected. Let $x_1,x_2,x_3$ and $x_4$ be 4 consecutive vertices on
$C_s^2$.  By Claim 3, $N_{G_i}(x_j)=V_i$ for $j=1,2,3,4$ and $i=1,2$.
Hence $C^2=x_1x_2P_1^2x_3x_4P_2^2C_s^2x_1$ is an $H^2$
of $G$.

%we are done by induction and using the P-C-C concatenation twice
%for each hamiltonian path square obtained in $G_i$\,(i=1,2).

So, we assume that $G[S]$ is 4-connected and $G[S]$ is a graph in
some of the exception families. In this case, we first show that
every graph in the exception families has a hamiltonian path square.
Then by concatenating the path square, $P_1^2$, and
$P_2^2$ together, we can get an $H^2$ of $G$.

Let $Q$ be a graph isomorphic to $(K_1\sqcup K_3)+(K_m\sqcup K_q)$
for some $m+q\ge 4$. We may assume, without loss of generality, that $m\ge 2$.
Then we let $P_3^2$ be a path square of $K_3$, $P_m^2$  a path
square of $K_m$, and $P_q^2$  a path square of $K_q$. Also, let $x$ be the single
vertex from $K_1$. Then $P_q^2P_3^2P_m^2x$ is a hamiltonian path square of $Q$.
The constructions for a hamiltonian path square for graphs in the families
of $(K_2\sqcup K_2)+(K_1\sqcup K_m)$,   $(K_2\sqcup K_3)+(K_1\sqcup K_m)$,
and  $(K_3\sqcup K_3)+(K_1\sqcup K_m)$ are similar, so we omit the details here.

Now let $P^2_s$ be a hamiltonian path square of $G[S]$,
and let  $x_1,x_2,x_3$ and $x_4$ be 4 consecutive vertices on
$P^2_s$.  By Claim~3, for any $v\in S$, $N_{G_i}(v)=V_i\,(i=1,2)$.
So $C^2=x_1x_2P_1^2x_3x_4P_2^2P_s^2x_1$ is an $H^2$
of $G$.

%Now, by using  P-P-C concatenation, we can first get an $H^2$ of $G[S\sqcup V_1]$.
%In the end, the P-C-C  concatenation can be applied on the
%$H^2$ of $G[S\sqcup V_1]$ and a hamiltonian path square $P_2^2$ of $G_2$.
%As $P_1^2$ has at least 4 vertices, we can insert $P_2^2$ such
%that the condition of P-C-C  is satisfied.

The remaining proof is divided into two cases according to the
connectivity of $G[S]$.  Let $G'=G[S]$.

\textbf{Case 1.}   Suppose $G'$ is connected but not 4-connected.

If $G'\cong K_4$, let $C_s^2=x_1x_2x_3x_4$ be an $H^2$ of it.
Then $C^2=x_1x_2P_1^2x_3x_4P_2^2x_1$ is an $H^2$
of $G$. So suppose $G'\not \cong K_4$.
As $|V(G')|\ge 4$ and $G'$ is not 4-connected,
$G'$ is not complete. Let  $S'$ be a minimum vertex-cut of $G'$.
Notice that $1\le |S'|\le 3$. Similar discussion
as in Claim~1-Claim~4 shows that $G'-S'$ has exactly two components, say,  $G_1'$ and
$G_2'$ such that each is a complete subgraph, and  $G'=G'[S']+(G'_1\sqcup G'_2)$.
 As $G'$ is also claw-free,
we see that $S'$ is $\overline{K_3}$-free.
Let $P_{1i}^2$ be a hamiltonian path square of $G_i'\,(i=1,2)$. Suppose,
without loss of generality, that $|V(P_{11}^2)|\le |V(P_{12}^2)|$.
We define two new vertex disjoint path squares of $G'$.

\begin{itemize}
  \item [C1.]$|S'|=1$. Let $S'=\{x_1\}$ and  $P_{21}^2=P_{11}^2x_1,\, P_{22}^2=P_{12}^2$;
  \item [C2.]$|S'|=2$.  Let $S'=\{x_1,x_2\}$ and  $P_{21}^2=P_{11}^2x_1,\,P_{22}^2=P_{12}^2x_2$;
  \item [C3.]$|S'|=3$.  Let $S'=\{x_1,x_2,x_3\}$,  and
assume that $x_1x_3\in E(G')$ by the fact that $S'$ is $\overline{K_3}$-free,
then let $P_{21}^2=x_1x_3P_{11}^2,\,P_{22}^2=P_{12}^2x_2$.
\end{itemize}

If C1 is true, then $\max\{|V_1|,|V_2|\}\ge 2$. Otherwise,
$S'\cup V_1\cup V_2$, a 3-set, separates $G_1'$ and $G_2'$,
contradicting the 4-connectedness assumption of $G$.
Assume, without  loss of generality, that $|V_1|\ge2$.
To specify the end vertices, we denote $P_{21}^2=x_1P_{21}^2x$ and
$P_{22}^2=zP_{22}^2x_2$, where $x_1\in S'$ and $z\in V(G')-S'$.
Clearly,  $x_1z\in E(G)$.
As $|V(G')|\ge 4$ and $|S'|=1$, $|V(P_{12}^2)|\ge 2$
by the assumption that $|V(P_{11}^2)|\le |V(P_{12}^2)|$.
Hence, both $P_{21}^2$ and $P_{22}^2$ have at least 2 vertices.
In specifying  one end of the hamiltonian path square $P_2^2$ of
$G_2$, let $P_2^2=P_2^2w$.
Then $x_1P_{21}^2xP_1^2x_2P_{22}^2zP_2^2wx_1$ is an $H^2$
of $G$ even if  $|V(P_{2}^2)|=1$.

For cases C2 and C3,
to specify the end vertices, we denote $P_{21}^2=x_1P_{21}^2x$ and
$P_{22}^2=zP_{22}^2x_2$, where $x_1,x_2\in S'$ and $x,z\in V(G')-S'$.
Since each of $P_{11}^2$ and $P_{12}^2$ has at least one vertex,
each of the $P_{21}^2$ and $P_{22}^2$ defined in C2 and C3 has
at least two vertices.
By the fact that $G'=G'[S']+(G'_1\sqcup G'_2)$ and the assumption that
$x_1x_3\in E(G)$, we see both
$P_{21}^2$ and $P_{22}^2$  are path squares satisfying $x_1z, xx_2\in E(G)$.
In specifying  one end of the hamiltonian path square $P_2^2$ of
$G_2$, let $P_2^2=P_2^2w$.
Then $x_1P_{21}^2xP_1^2x_2P_{22}^2zP_2^2wx_1$ is an $H^2$
of $G$ even if $|V(P_{1}^2)|=1$ or $|V(P_{2}^2)|=1$.

\textbf{Case 2.}   Suppose $G'$ is disconnected.

As $G$ is claw-free and $G=G'+(G_1\sqcup G_2)$, we see that $G'$
consists of exactly two complete components, say $G_1'$ and $G_2'$.
So $G=(G_1'\sqcup G_2')+(G_1\sqcup G_2)$ and $V_1\cup V_2$
is also a vertex-cut of $G$.  For $i=1,2$, let $|V(G_i')|=|V_i'|$.
So $|V_1\cup V_2|\ge |V_1'\cup V_2'|=|S|$, by the minimality of $|S|$.
Recall that $G$
is not isomorphic to any of the graphs in the following families:
\begin{enumerate}
  \item [$($i$)$] $(K_1\sqcup K_3)+(K_m\sqcup K_q)$ with $m+q\ge 4$;
  \item [$($ii$)$] $(K_2\sqcup K_2)+(K_1\sqcup K_m)$ with $m\ge 3$;
  \item [$($iii$)$] $(K_2\sqcup K_3)+(K_1\sqcup K_m)$ with $m\ge 3$;
  \item [$($iv$)$]$(K_3\sqcup K_3)+(K_1\sqcup K_m)$ with $m\ge 3$.
 \end{enumerate}

 Assume first that $\min\{|V_1|, |V_2|, |V_1'|, |V_2'|\}\ge 2$.
 In specifying  the end vertices, we let $P_1^2=x_1P_1^2y_1$,
 $P_2^2=x_2P_2^2y_2$, $P_{11}^2=x_{11}P_{11}^2y_{11}$, and
 $P_{12}^2=x_{21}P_{12}^2y_{21}$
 be the hamiltonian path square of $G_1$, $G_2$, $G_1'$ and $G_2'$,
 respectively.
 Then as $G=(G_1'\sqcup G_2')+(G_1\sqcup G_2)$, we know
 $x_1P_1^2y_1x_{11}P_{11}^2y_{11}x_2P_2^2y_2x_{21}P_{12}^2y_{21}x_1$
 is an $H^2$ of $G$.
 So assume, without loss of generality,
 that $|V_1|=1$.
 Then as $G$ is not isomorphic to any graphs in (i)
 and $|V_1\cup V_2|\ge |V_1'\cup V_2'|=|S|\ge 4$, we have that $|V_2|\ge 4$.
 So, $G_1\sqcup G_2\cong K_1\sqcup K_m$ for some $m\ge 4$.
 Also, as $G$ is not isomorphic to any graphs in (i)-(iv), $G_1'\sqcup G_2'\not \cong
 K_1\sqcup K_3, K_2\sqcup K_2, K_2\sqcup K_3, K_3\sqcup K_3$.
 This indicates that $\max\{|V_1'|, |V_2'|\}\ge 4$.  We may assume,
 without loss of generality, that $|V_1'|\ge 4$.
 Let $P_2^2=x_{21}x_{22}\cdots x_{2,s-1}x_{2s}\,(s\ge 4)$
 be the hamiltonian path square of $G_2$ specified earlier,
 $P_{11}^2=x_{11}x_{12}\cdots x_{1,t-1}x_{1t}\,(t\ge 4)$
be a hamiltonian path square of $G_1'$,  and let
$P_{12}^2$ be a hamiltonian path square of
$G_2'$.
%Then $$
%
%We have that either (a) $\min\{|V_1|,|V_2|\}\ge 2$ and $\min\{|V(G_1')|,|V(G_2')|\}\ge 2$
%or (b) $|V(G_i')|=1, |V(G_{3-i}')|\ge 4$ and $|V_1|+|V_2|\ge |V(G_i')|+|V(G_{3-i}')|$ for
%$i=1,2$. (Notice that $|V_1|+|V_2|\ge |V(G_i')|+|V(G_{3-i}')|$ follows from the minimality of $S$.)
%If (a) is the case, then we can find an $H^2$ of $G$ by Remark 3. So we assume (b)
%is the case, and assume without loss of generality that
%$|V(G_1')|\ge 4, |V(G_{2}')|=1$ and $|V_1|\ge |V_2|$. Furthermore, as
%$G$ is not isomorphic to any graphs in $($ii$)$, $($iii$)$, and $($iv$)$,
%we have $|V_1|\ge 4$. Let  $P_1^2=x_{11}x_{12}\cdots x_{1,t-1}x_t\,(t\ge 4)$
%be a hamiltonian path square of $G_1'$, $P_2^2$ be a hamiltonian path square of
%$G_2'$, $P_3^2=x_{31}x_{32}\cdots x_{3,s-1}x_s\,(s\ge 4)$
%be a hamiltonian path square of $G_1$,  and $P_4^2$ be a hamiltonian path square of
%$G_2$.
Then $x_{11}x_{12}P_1^2x_{13}x_{14}P_{11}^2x_{1,t-1}x_{1t}x_{21}x_{22}P_{12}^2x_{23}x_{24}P_2^2 x_{2,s-1}x_{2s}x_{11}$
is an $H^2$ of $G$.

The proof of Theorem~\ref{clawp4} is then complete.
\qqed

{\noindent \bf Acknowledgements:}
The authors wish to thank  the two anonymous referees for their helpful comments.

\bibliographystyle{plain}
\bibliography{SSL-BIB}

\end{document}